# A Note on a Noetherian Fully Bounded Ring .

**By C.L.Wangneo**

**Jammu ( J& K)**

**India -180002 .**

**Abstract :** We prove the following theorem ;

**Theorem :--Let R be a_ a prime noetherian ring with ]R] = n , n a finite nonnegative integer . For a fixed integer m, 0 ≤ m < n call Xm = { All p ε spec. R / k-dim. (R/ P) = m } , the full set of m-prime ideals of R. Let cm ={c ε R / k-dim (R./cR) < m } and let c= {Right ideals I ≤ R / I п cm ≠ ϕ } .Let g = {Right ideals J ≤ R / k-dim . (R/J) < m} be the m-gabriel filter and let h = { Right ideals K ≤ R / k п c(p) ≠ ϕ , for each p ε xm } .Let vm = п c(p) , for all p ε xm and let v== { Right ideals I ≤ R / I п vm ≠ ϕ }. We say xm has the right intersection condition if h=v. . Suppose any m critical right R module M with Ass. M = p is such that IMI = IR/pI=m (or M is a torsionfree right R/P –module ) .Then the following conditions are equivalent ;**

**(a) xm has the right intersection condition .**

**(b) (i) g =v . In this case then vm is a right ore set in R . Let Rv denote the quotient ring of R at the right ore set vm .**

**(ii) Moreover then any simple Rv module namely RV/Jv , Jv a maximal right ideal of Rv with r(Rv/Jv) = qv is a torsion free Rv / qv module.**

**(c) (i) g=c . Then cm is a right ore set in R . Let Rc denote the quotient ring of R at the right ore set cm .**

**(ii) Moreover then any simple Rc module namely Rc/Jc , Jc a maximal right ideal of Rc with r(Rc/Jc) = qc is a torsion free Rc / qc module**

**The theorem is proved under the( underlined) weaker hypothesis on the prime noetherian ring R than for a prime noetherian ring that is either fully bounded or has the bijective Gabriel correspondence . For these classes of rings the theorem remains true always for all m , 0≤ m≤ n . Moreover the theorem is true if we replace k-dim. R =n , n finite by any ordinal number .**

**Introduction** :- This paper is divided into two sections . In section (1) we introduce some preliminary definitions and results for a prime , noetherian ring R of krull dimension , n , where $n \geq o$ is a finite integer and which will be of use to us throughout this paper. Then in section (2) we characterise the intersection condition of a full set of prime ideals xm of a noetherian prime ring R . We sum up all this in the form of our main theorem at the end of section (2) which states the following ;

Main Theorem :-- Let R be a prime noetherian ring with ]R] = n , n a finite nonnegative integer . For a fixed integer m, $0 \leq m < n$ call Xm = { All p ε spec. R / k-dim. (R/ P) = m } , the full set of m-prime ideals of R. Let cm ={c ε R / k-dim (R./cR)I < m } and let c= { I ≤ R / I п cm ≠ ф } .Let g = { I ≤ R / k-dim . (R/I) < m} be the m-gabriel filter and let h = { J ≤ R / J п c(p) ≠ ф , for each p ε xm } .Let vm = п c(p) , for all p ε xm and let v== { J ≤ R / J п vm ≠ ф ,} We say xm has the right intersection condition if h=v. . Suppose any m critical right R module M with Ass. M = p is such that IMI = IR/pI=m (or M is a torsionfree right R/P –module ) .Then the following conditions are equivalent ;

(a) xm has the right intersection condition .

(b) (i) g =v . In this case then vm is a right ore set in R . Let Rv denote the quotient ring of R at the right ore set vm .

(ii) Moreover then any simple Rv module namely RV/Jv , Jv a maximal right ideal of Rv with r(Rv/Jv) = qv is a torsion free Rv / qv module.

(c) (i) g=c . Then cm is a right ore set in R . Let Rc denote the quotient ring of R at the right ore set cm .

(ii) Moreover then any simple Rc module namely Rc/Jc , Jc a maximal right ideal of Rc with r(Rc/Jc) = qc is a torsion free Rc / qc module

**The theorem is proved under the( underlined) weaker hypothesis on the prime noetherian ring R than for a prime noetherian ring that is either fully bounded or has the bijective Gabriel correspondence . For these classes of rings the theorem remains true always for all m , 0≤ m≤ n . Moreover the theorem is true if we replace k-dim. R =n , n finite by any ordinal number .**

## Notation and Terminology :

**(a) Throughout in this paper a ring is meant to be an associative ring with identity which is not necessarily a commutative ring.**

**We briefly mention below what few definitions and terms occur in this paper . For all these definitions and terms and the preliminary results one may consult [1] and [2] .**

**(b) By a module M over a ring R we mean that M is a right R-module unless stated otherwise. For the basic definitions regarding right noetherian modules over right noetherian rings and the definitions of uniform modules , right Krull dimension , critical modules , krull homogenous modules we refer the reader to both [1] and [2 ]. Also one may consult both [1] and [2[ for the definition and properties of Ass.(M) and its associated concepts such as faithful and fully faithful module over a right noetherian ring .We will throughout call a finitely generated right module M a m-critical module if (i) k-dim. (M) = m and (ii) M is a critical right module. Moreover we will use the following terminology throughout .**

**(c) If R is a ring then we denote by Spec.R , the set of prime ideals of R . Also right Krull dimension of the right R-moodule M if it exists will throughout be denoted by |M| . For two subsets A and B of a given set , A ≤ B means B contains A and A < B denotes A ≤ B but A≠B. Also for two sets A and B , A⊄B denotes the set B that does not contain the subset A. For an ideal A of R , c(A) denotes the set of elements of R that are regular modulo A . Recall if S ≤ M is a non-empty subset , then we denote the right**

annihilator of S in R by r (S) . Also recall that if R is a ring and M is a right R module and if T is a multiplicatively closed subset of regular elements of R then a submodule N of M is said to be T tosion if for any element x of M there exists an element t in T such that $x t = 0$ . M is said to be T torsion free if for any nonzero element x of M and for any t in T , $x t \neq 0$ . If M is not T- torsion , then there exists at least one non-zero element x of M such that x is a T Torsion free element of M . This means $x t \neq 0$ for some element $t \neq 0$ of T . Similarly if M is not T-torsionfree , then there exists at least one non-zero element y of M such that y is a T- Torsion element of M . This means $yt = 0$ , for at least one element $t \neq 0$ in T . In this case then we say that M is a T-non-torsion right R module .

(d) We will also use the following notation throughout ;

(i) Spec.(R) = family of prime ideals of R .

(ii)max. Spec.(R) = family of maximal ideals of R .

(iii) prim Spec.(R) = family of right primitive ideals of R .

Note that for a ring R max.spec. (R) ≤ prim .spec. (R) .

We briefly mention that if S is an extension ring of the ring R we denote by r(W) the right annihilator in S of any simple S- module (which is a prime ideal of S ).

A word about the numbering of a theorem , propostion ,lemma etc. All such numbering occurs on the left of each proposition , lemma etc. but while referring to these theorems lemmas etc. we refer them in the usual way in the paper For example (1.2) Lemma would be referred as Lemma (1.2) in the paper while quoting it .in the paper .

(1) section ( preliminaries ) :- Let R be a prime noetherian ring with IRI = n, n a nonnegative integer.We now state some preliminaries regarding this paper . All these preliminaries can be found in [2]

**(1.1)Definition and Notation:** (a) Let R be a prime , noetherian ring . Let ]R] = n , n , a finite integer , n ≥ 0 . For a fixed integer m , 0 ≤ m < n , let xm = { All p ε spec. R / |R/ P| = m } We call xm as the full set of m-prime ideals of R .Then we let vm=∩ c(p) ; for each p ε xm. and let cm = { all c ε R / IR/cR I< m } Then cm is a multiplicatively closed subset of regular elements of R and cm ≤ vm and vm is also a multiplicatively closed subset of regular elements of R

(b)For the prime noetherian ring R with IR I= n , n ≥ 0 let xm , vm and cm be defined as in definition (a) above . We define the following families of right ideals of R , namely ;

c = { Right ideals I of R / I∩Cm ≠ϕ} .

v = { Right ideals I of R / I∩ Vm ≠ϕ} .

g = { Right ideals I of R / IR/I I< m }.

h = { Right ideals I of R / I∩ c(p) ≠ϕ for each pεXm}

(c) For the prime noetherian ring R with IR I=n , n ≥ 0 if xm , vm,cm and v , c , g and h are as in definitions (a) and (b)above then We say xm has the right intersection condition if for any right ideal I of R , I∩ c(p) ≠ϕ for each pε Xm, implies that I∩ Vm ≠ϕ. Note this is equivalent to saying that V=h, because v ≤h always

**(1.2) Lemma :-** Let R be a prime,noetherian ring with ]R] =n, n, a finite integer , n ≥0.For a fixed integer m, 0 ≤ m < n let xm, vm , cm be as in the above definition (1.1) (a) . Also Consider the families of right ideals c, v, g and h as in definition(1-1)(b) . Then the following are true ;

(a) cm ≤ vm and both cm and vm are multiplicatively closed subsets of R

**(b) $C \leq g \leq h$ . Also $c \leq v \leq h$.**

**(c) If Xm has the right intersection condition then we have**

**(i) $v = h$ and**

**(ii) $g \leq v$ .**

**(1.3) Definition (Extension ring of a ring R ) :- Let R be a prime noetherin ring with goldie quotient ring Q and let S be an extension ring S of R contained in Q such that $1_R = 1_s = 1$. Throughout we will write $S = R^e$ Then we define the following ;**

**(a) Definition :- (i) For any right ideal A of R , $A^e$ is the right ideal of S such that $A^e = AS$ , that is $A^e$ is the right ideal generated by A in S . $A^e = AS$ . We call $A^e$ the extension of the right ideal A to $S = R^e$**

**(b) For any right ideal B of S , $B^C = B \cap R$ . Clearly $B^C$ is a right ideal of R . We call $B^c$ the contraction of the right ideal B of $S=R^e$ to R**

**(1.4 ) Gabriel filters :- We now connect ore-sets with Gabriel filters. -**

**(a) Definition :- Let R be a ring . Let g = { I < R / I is a right ideal } . Let g be a non-empty family . Then g is called a m-Gabriel filter if g has the following properties ;**

**(i) For any I ε g if J is a right ideal such that $I \leq J$ , then J Ɛ g .**

**(ii) For I, J in g , $I \cap J$ Ɛ g .**

**(iii) For any a in R , and I ε g , the right ideal a-1 I = { x in R / ax in I } , is also in g .**

**Note if J = a-1 I , then J is a maximal right ideal of R such that a J < I .**

(b) Definition_ In a ring R a multiplicative Gabriel filter is a nonempty family g of right ideals of R such that

(i) g is a Gabriel filter

(ii) For I, J in g the product of the right ideals IJ Ɛ g .

(1.5) Remark :- If cm = {c in R / IR./ cRI < m } and c= { all right ideals k in R / K п cm ≠ ф } then g= c implies in particular that cm is a right oreset of R . However given cm as above let w= { all right ideals J in R / R/J is a right cm torsion R- module } then w is a Gabriel filter and w ≤ c where c= { all right ideals k in R / k п cm ≠ ф } Also it is not difficult to see that cm is a right ore set if and only if w=k . Clearly g=c is a stronger condition than w=k . Therefore we give the following definition first :-

(1.6)Definition :- Let cm , c, g and w be as in Remark (2.2) above . Then we have the following definition ;

We say cm is a right ore set if w=c and cm is a strongly right ore set if g =c where g is an m-gabriel filter as defined above .. It is true that g=c implies w=c always .

Section (2) (The intersection condition ) :- In this section we prove our main theorem for rings that are far more general than fully bounded noetherian rings or noetherian rings that have the bijective Gabriel correspondence . So it would be appropriate to say a few words about these rings to conform to the title of this paper .

(2.1) Remark :- For the definition of a fully bounded noetherian ring and a noetherian ring with the bijective Gabriel correspondence we refer the reader to chapter (8) 0f [2]( p.132 ) and theorem 8.14 therein as well as to chapter (8) 0f [1] (p.51) and theorem (8.6) therein . Moreover we need to mention that in this paper we will be using only the following property of the rings mentioned above namely if R is a prime noetherian ring with k-dim. (R) =n , $n \geq 0$ then for any fixed non-negative integer m , $0 \leq m \leq n$ we have that every m- critical right R module M with Ass. M = p is such that IMI = IR/pI=m (or equivalently M is a torsionfree R / P –module ) . We must also mention that we will work henceforth throughout this paper with this hypothesis underlined always .

We will now prove our main theorem namely theorem (2.4) . Towards this end we state and prove the following results below

(2.2) proposition :- Let R be a prime , noetherian ring with ]R] = n , n ( n a finite nonnegative integer ). For a fixed integer m, $0 \leq m < n$ let Xm = { All p ε spec. R /|R/ P| = m } , be the full set of m-prime ideals of R. Let cm ={c ε R / IR./cRI $<$ m } .Let g = { I ≤ R / IR/II $<$ m} be the m-gabriel filter and let h = { J ≤ R / J п c(p) ≠ ϕ } .Let vm = {п c(p) , for all p ε xm }. Then the following are true ;

(a) (i) c ≤ g ≤ h .

(ii) c ≤ v ≤ h .

(b) If xm has the right intersection condition then we have ;

(i) v=h and g ≤ v .Hence c≤ g ≤ v =h .

(2.3)Remark :- In proposition (2.2) above we consider two cases namely ;

Case (i) :- Assume g $<$ h ( strict containment ). In this case we choose a right ideal J ε h maximal such that J is not in g . Then R/ J is a m- critical right R – module such that for any factor module of R/J by a right ideal I properly containing J we have IR/II $<$ m , hence I ε g ..

Case (ii) : -- Assume c $<$ g (strict containment ) .In this case we choose a right ideal J ε g maximal such that J is not in c . Then R/ J is a k- critical right R – module with k $<>$ m such that for any factor module of R/J by a right ideal l properly containing J we have l п cm ≠ф .

In case (ii) above if we consider the extension $J^e$ of the right ideal J to an extension ring $R^e$ of R whose set of units is the set cm then clearly $J^e$; is a maximal right ideal of $R^e$ and hence $R^e$ /$J^e$ is a simple right $R^e$ – module . However the same cannot be said regarding $R^e$ /$J^e$ in case (i) above .

Intersection condition means given c, g, v, and h we have to prove that v=h . Under certain conditions it would be a help if we know that g=h . The theorem below proves this under a condition which we have underlined.

**(2.4) Theorem :-** Let R be a_ a prime noetherian ring with ]R] = n , n a finite nonnegative integer . For a fixed integer m, 0 ≤ m < n let Xm = { All p ε spec. R /|R/ P| = m } , be the full set of m-prime ideals of R. Let cm ={c ε R / IR./cRI < m } .Let g = { I ≤ R / IR/II < m} be the m-gabriel filter and let h = { J ≤ R / J п c(p) ≠ ф } .Let vm =п c(p) , for all p ε xm . If every m- critical right R module M with Ass. M = p is such that IMI = IR/pI=m (or M is a torsionfree R / P – module ) then we have that g=h .

**Proof :--** First note from lemma (2.2) we have that c ≤ g ≤ h . Let us suppose that g ≠ h . Then we have that g < h (strict containment ) . . Choose J in h such that J is not in g . We may choose J to be max. such . Then R/J is an m-critical right R- module because any factor module of R/J by a right ideal l of R properly containing J is such that I ε g and hence IR/I I <m . Therefore by the underlined hypothesis of the theorem IR/JI = IR/pI . Hence p εxm .But J is in h and hence J п c(p) ≠ ф, for each p εxm by the definition of h . Thus IR/JI < m showing thereby that J ε g . Hence g = h .

**(2.5) Main Theorem :--** Let R be a_ a prime noetherian ring with ]R] = n , n a finite nonnegative integer . For a fixed integer m, 0 ≤ m < n call Xm = { All p ε spec. R / k-dim. (R/ P) = m } , the full set of m-prime ideals of R. Let cm ={c ε R / k-dim (R./cR)I < m } and let c= { I ≤ R / I π cm ≠ ϕ } .Let g = { I ≤ R / k-dim . (R/I) < m} be the m-gabriel filter and let h = { J ≤ R / J π c(p) ≠ ϕ , for each p ε xm } .Let vm = π c(p) , for all p ε xm and let v== { J ≤ R / J π vm ≠ ϕ ,} . We say xm has the right intersection condition if h=v. . Suppose any m critical right R module M with Ass. M = p is such that IMI = IR/pI=m (or M is a torsionfree right R/P –module ) .Then the following conditions are equivalent ;

(a) xm has the right intersection condition .

(b) (i) g =v . In this case then vm is a right ore set in R . Let Rv denote the quotient ring of R at the right ore set vm .

(ii) Moreover then any simple Rv module namely RV/Jv , Jv a maximal right ideal of Rv with r(Rv/Jv) = qv is a torsion free Rv / qv module.

(c) (i) g=c In this case then cm is a right ore set in R . Let Rc denote the quotient ring of R at the right ore set cm .

(ii) Then any simple Rc module namely Rc/Jc , Jc a maximal right ideal of Rc with r(Rc/Jc) = qc is a torsion free Rc / qc module

Proof :- (a) => (b)(i) Note first that using Theorem (2.4) above we have g=h and since v ≤ h , hence we must have that v ≤ g . From prlemma (2.2) above it is clear that if xm has the right intersection condition then g ≤ v . therefore from above we get that g=v .

(a) => (b)(ii) We now prove (b) (ii) . So let Jv be a maximal right ideal of Rv so that the module Rv/Jv is a simple Rv module .We show that Rv/Jv is a torsionfree Rv /qv module where qv = r(Rv/Jv) . To see this observe that since Jv is a proper right ideal of Rv , hence J π v = ϕ and since Jv is a maximal right ideal of Rv thus J is maximal such that J π v = ϕ . Thus J does not belong to v . But g=v . Hence J does not belong to g . Thus IR/JI ≥ m . Moreover any right ideal k containing J properly is such that k π v ≠ ϕ ( use the fact that

**Jv is a maximal right ideal of Rv ) and g=v implies k ε g . Thus IR/kI < m . Thus IR/JI is a m-crtical right R-module. Hence by the underlined hypothesis on the ring R we have that R/J is a torsion free R/q module where q=Ass. (R/J ) . Hence IR/JI = I R/q I =m (in fact q ε xm )**

**(b) => (c) Given (b) , namely g=v , we see that then g=c and vm=cm and so v=c .Therefore to see the truth of (c) all that we have to do is to replace g=v by g=c .**

**(c) => (a): Let Rc be the ring R localised at the right ore set cm . Clearly the set cm is the set of units of Rc . To prove the right intersection condition for xm we prove v=h or v=h=c . Suppose not . Then let v<h . But g=c=v and v< h implies g< h (strict containment) . Choose a right ideal J in h such that J is not in v and J is maximal such . . Hence J is not in c=g . Thus R /J is an m-critical right R-module. Clearly by the undelined hypothesis of the theorem if P = Ass. R/J , then R/J is a torsion free R/p- module . Hence p ε xm . Since J is in h hence J π c(p) ≠ ф . Thus IR/JI < m implies thereby that J is in g and g=c . Hence we must have that v = h which proves the right intersection condition .**

**(2.6) Theorem :- Let R be a prime noetherian ring satisfying all the conditions of theorem (2.5) . Then the conclusion of the above theorem remain true for all m , 0≤ m≤ n , if in particular R is either a noetherian fully bounded ring or a noetherian ring having the bijective Gabriel correspondence .**

**Proof :- If the ring R is either a noetherian fully bounded ring or a noetherian ring having the bijective Gabriel correspondence then the underlined hypothesis of theorem (2.5) above holds true for a m-critical R- module W of R for any m . Hence in particular the theorem is true for all such m , 0≤m ≤ n**

**(2.7) corollary :--_ In theorem (2.5) above If xm has the right intersection condition then we have the following further as true ;**

**(a) max. spec. (Rv ) = {Pv / for all P ε x } .= Prim. Spec(Rv) and for all p ε x IRv / pv I = 0 .**

**(b) max. spec. (Rc ) = {Pc / for all P ε x } .= Prim. Spec(Rc) and for all p ε x IRc / pc I = 0 .**

**(2.8)Remark :- It is now appropriate to refer to the definition in [2] of exercise 12 of chapter 12 named as localizability of x. In fact it is used there more generally as proved in [2] exercise 12 of chapter 12 to characterise the right intersection condition of x which is an arbitrary set of incomparable prime ideals of a noetherian ring .**